\documentclass[preprint,12pt]{elsarticle}

\usepackage{amssymb}
\usepackage{amsthm}
\usepackage{amsmath}

\newtheorem{remark}{Remark}
\newtheorem{theorem}{Theorem}

\newcommand{\bB}{\pmb{B}}
\newcommand{\bV}{\pmb{V}}
\newcommand{\bv}{\pmb{v}}
\newcommand{\bx}{\pmb{x}}
\newcommand{\bA}{\pmb{A}}
\newcommand{\bC}{\pmb{C}}
\newcommand{\bF}{\pmb{F}}
\newcommand{\bU}{\pmb{U}}
\newcommand{\bP}{\pmb{P}}
\newcommand{\bsigma}{\pmb{\sigma}}
\newcommand{\bff}{\pmb{f}}
\newcommand{\bu}{\pmb{u}}
\newcommand{\iT}{\mathcal{T}}

\newcommand{\iK}{\mathcal{K}}
\newcommand{\iN}{\mathcal{N}}

\newcommand{\mR}{\mathbb{R}}
\newcommand{\Div}{\text{div}\;}
\begin{document}

\begin{frontmatter}



\title{A staggered finite element cell-centered scheme for nearly incompressible elasticity on general meshes}


\author[Phuong]{T.T.P. Hoang}
\address[Phuong]{Department of Mathematics, Ho Chi Minh city University of Pedagogy, 280 An Duong Vuong street, ward 4, district 5, Ho Chi Minh city, Vietnam \vspace{5pt}}
\author[aHai]{Ong Thanh Hai}
\address[aHai]{Faculty of Mathematics and Computer Science, University of Science, VNU HCMC, 227 Nguyen Van Cu Street, District 5, Ho Chi Minh City, Vietnam.}  
\author[aHai]{ H. Nguyen-Xuan}
\begin{abstract}
We apply the finite element cell-centered (FECC) scheme \cite{PO12} to the solution of the nearly incompressible elasticity problem. By applying a technique of dual mesh, such a low-order finite element scheme can be constructed from any given mesh and it is proved to satisfy the ``macroelement condition" \cite{RS84}, meaning that the stability condition is fulfilled. Numerical results show that the method, which is simple to implement, is effective in terms of accuracy and computational cost compared with other methods. 
\end{abstract} 

\begin{keyword} nearly incompressible elasticity \sep finite elements \sep cell-centered scheme \sep macroelement condition \sep inf-sup condition



\end{keyword}

\end{frontmatter}


%
%
%
%
%
%
%
%
%
%
%
%
\section{Introduction}
\label{Sec:intro}
In this paper, we present a new numerical method for nearly incompressible elasticity problems using low-order finite elements. 
The scheme, firstly introduced for solving stationary diffusion problems, has many advantages: i) it can deal with general meshes and it involves only cell unknowns; ii) it gives higher accuracy while the computational cost is the same as other cell-centered scheme; iii) the stability condition is obtained by using the ``macroelement condition" instead of using bubble functions (as in \cite{BL05}), which means that the scheme is simple to implement and very well-suited to the problem. 

The rest of this paper is organized as follows: in Section~\ref{Sec:model}, a two-dimensional model problem in mixed form 
is introduced. In Section~\ref{Sec:discrete}, we present the discretizations using FECC scheme (in which a technique of dual mesh is employed) and formulate the corresponding discrete problem. Some theoretical results concerning the stability of the new scheme are proved, and the advantages of using FECC scheme in terms of accuracy and computational cost are also discussed. Numerical experiments for two-dimensional problems comparing the performance of different schemes are shown in Section~\ref{Sec:NumRe}.  

%
%
\section{A model problem}
\label{Sec:model}
For a bounded domain $ \Omega $ in $ \mR^2 $ with Lipschitz boundary $ \partial \Omega $, we consider the following stationary linear elasticity equation:
\begin{equation} \label{primal}
-\Div \bsigma(\bu)=\bff \quad \text{in} \; \Omega,
\end{equation}
where $ \bu $ is the displacement of an elastic material, $ \bsigma $ the Cauchy stress and $ \bff $ the body forces. For simplicity, we impose a homogeneous Dirichlet boundary condition on $ \partial \Omega$:
\begin{equation} \label{BC}
\bu = 0 \quad \text{on} \; \partial \Omega.
\end{equation}
If the elastic material is isotropic and linear, the stress $ \bsigma (\bu) $ is defined by
\begin{equation*}
\bsigma (\bu)=2\mu \varepsilon (\bu) + \lambda \Div \bu \mathbf{Id},
\end{equation*}
where $ \mathbf{Id} $ is the identity matrix of size $ 2 $, $ \varepsilon (\bu) $ the infinitesimal strain tensor defined by
\begin{equation*}
\varepsilon (\bu) = \frac{1}{2}\left (\nabla \bu + (\nabla \bu)^{T}\right ),
\end{equation*}
(here $ A^{T} $ denotes the transpose of a matrix A),  $ \lambda $ and $ \mu $ are the Lam\'e constants defined by
\begin{equation*}
\lambda =\frac{\nu E}{(1+\nu)(1-2\nu)}, \quad \mu=\frac{E}{2(1+\nu)},
\end{equation*} 
with $ \nu $ the Poisson's ratio and $ E $ the Young's modulus. In this work, we study the case where $ \nu $ is close to $ 0.5 $ (or $\lambda$ is large), i.e. the material is nearly incompressible (for example rubber or rubber-like materials). It is well-known that for such a case, the standard finite elements might give inaccurate results due to volumetric locking and instability. In order to avoid this, one may use the mixed formulation by introducing an additional variable, the pressure $ p $, as follows:
\begin{equation*}
p:=\lambda \Div \bu \quad \text{in} \; \Omega,
\end{equation*} 
 and rewrite problem \eqref{primal}-\eqref{BC} equivalently in a mixed displacement-pressure form as
  \begin{subequations} \label{mixed}
 \begin{align}
 -\Div \left (2\mu \varepsilon (\bu) + p \mathbf{Id}\right ) &= \bff & \text{in} \; \Omega,  \label{1stmixed} \\
 \Div \bu - \frac{1}{\lambda} p & = 0 & \text{in} \; \Omega, \label{2ndmixed} \\
 \bu & =0 & \text{on} \; \partial \Omega.  \label{3rdmixed}
 \end{align}
\end{subequations}
To derive the variational formulation of \eqref{mixed}, we first introduce the following Sobolev spaces
\begin{equation*}
\bV_{0}=\left (H_{0}^{1}(\Omega)\right )^{2} \; \; \text{and} \; \; L_{0}^{2}(\Omega):=\left \{ q \in L^{2}(\Omega): \int_{\Omega} q d\Omega = 0\right \}.
\end{equation*}
Denote by $\| \cdot \|_{0}$ and $\| \cdot \|_{1}$ the norms defined on $L_{0}^{2}$ and $\bV_{0}$ respectively. We  shall seek for $ \bu \in \bV_{0} $ and $ p \in L_{0}^{2}(\Omega) $. The latter is obtained by integrating  equation \eqref{2ndmixed} over $ \Omega $ and using divergence theorem together with boundary condition \eqref{3rdmixed}. We also define the bilinear forms: 
\begin{equation*}
\begin{array}{lll}
a:& \bV_{0} \times \bV_{0} &\rightarrow \mR \\
& (\bu, \bv) & \mapsto a(\bu, \bv) = 2\mu \int_{\Omega} \varepsilon (\bu) : \varepsilon(\bv) \, d\Omega, \vspace{3pt}\\
b:& \bV_{0} \times L_{0}(\Omega) & \rightarrow \mR \\
& (\bu, q) & \mapsto b(\bu, q) = \int_{\Omega} q \Div \bv \, d\Omega, \vspace{3pt} \\
c:& L_{0}(\Omega) \times L_{0}(\Omega) & \rightarrow \mR \\
& (p, q) & \mapsto c(p, q) = \int_{\Omega}  pq \, d\Omega, \vspace{3pt}\\
L_{f}:&  \bV_{0} & \rightarrow \mR \\
& \bv & \mapsto L_{f}(\bv) = \int_{\Omega}  \bff \cdot \bv \, d\Omega,
\end{array}
\end{equation*}
where $ \bff \in \left (L^{2}(\Omega)\right )^{2} $. \vspace{4pt}

With these notations, the variational form of \eqref{mixed} is written as:  \vspace{-0.2cm}
\begin{eqnarray} \label{variational-mixed}
\mbox{Find $ \bu \in \bV_{0} $ and $ p \in L_{0}^{2}(\Omega) $ such that}
   \hspace{4cm}\nonumber\\    
\begin{array}{rll} a(\bu, \bv) + b(\bv, p)   & = L_{f}(\bv ), & \forall \bv \in \bV_{0}, \vspace{3pt}\\
 b(\bu, q) - \frac{1}{\lambda } c(p,q) & = 0, & \forall q \in L_{0}^{2}(\Omega).
\end{array}
\end{eqnarray}
In the next section, we present the finite element cell-centered (FECC) scheme \cite{PO12} to obtain a numerical solution to problem \eqref{variational-mixed}.
%
%
\section{The discrete problem using FECC}
\label{Sec:discrete}
We extend the idea of FECC scheme for stationary diffusion problems \citep{Othesis12} to the case of nearly incompressible elasticity problems. We recall that FECC scheme is based on a technique of dual mesh and it involves only cell unknowns. Three meshes are respectively constructed, the primal mesh, the dual mesh and the third mesh. The third mesh consisting of triangular elements can be seen as a refinement of the primal mesh. However, the number of unknowns does not increase compared with other cell-centered schemes. Thus, FECC scheme is effective in terms of accuracy and computational cost. Furthermore, by choosing appropriate approximate spaces for $\bu$ and $p$ we obtain a stable, low-order finite element scheme. This is done using macroelement techniques which turn out to be very natural in this case. 

In the following, we construct the meshes and introduce the approximate spaces to obtain the discrete problem associated with \eqref{variational-mixed}(cf. Subsection~\ref{Subsec:Discrete}). Then we prove that the discrete problem is well-posed and the scheme is stable (cf. Subsection~\ref{Subsec:proof}). Finally, we perform some calculations to obtain the algebra system corresponding to the discrete problem and discuss some issues concerning the efficiency of the scheme (cf. Subsection~\ref{Subsec:efficiency}).

\subsection{The discretization} \label{Subsec:Discrete}
For a polygonal domain $ \Omega $, consider a triangulation $ \iT_{h} $ (of $\Omega$) that consists of non-empty connected close disjoint subsets of $\Omega$: $$ \overline{\Omega} = \bigcup_{K \in \iT_{h}} K.
$$
We assume that each element $K \in \iT_{h}$ is a star-shaped polygon in which we will choose a point $C_{K} \in \text{int}(K)$ and call it the mesh point of $K$. Throughout the paper, we refer to $\iT_{h}$ as the primal mesh. Next, we briefly recall the construction of the dual mesh $\iT_{h}^{*}$ and the third mesh $\iT_{h}^{**}$ which will be necessary to define the FECC scheme (see \citep{PO12, Othesis12} for a more detailed presentation). To define the dual mesh, we assume that the line joining two mesh points of any two neighboring elements is inside $\Omega$ but it doesn't need to intersect the common edge of the two element since the problem is homogeneous. 

The introduction of the dual mesh $\iT_{h}^{*}$ is based on the primal mesh so that each dual control volume  of $\iT_{h}^{*}$ corresponds to a vertex of $\iT_{h}$. Denote by $\iN$ the set of all nodes or vertices of $\iT_{h}$
$$ \iN : = \left \{ i: \; \text{i is a vertex of element $K \in \iT_{h}$ } \right \}.
$$
For each $i \in \iN$, denote by 
$$\iT_{i}:=\left \{ K \in \iT_{h}: \; \text{K shares the vertex i} \right \},$$ 
the set of primal elements that have $i$ as their vertex. We consider two cases (see Figure): 
\begin{enumerate}
\item[(a)] If $i$ is an interior vertex, then by connecting mesh points of neighboring elements in $\iK_{i}$ we obtain the dual control volume $M_{i}$ associated with the vertex $i$.
\item[(b)] If $i$ is on the boundary $\partial \Omega$, denote by $E_{i} \subset \partial K_{i}^{E}$ and $F_{i} \subset \partial K_{i}^{F}$ the two edges on the boundary that have $i$ as their vertex in which $K_{i}^{E}, K_{i}^{F} \in \iK_{i}$ can be two distinguished elements or just one element. 
The dual control volume $M_{i}$ is defined by joining mesh points of neighboring elements in $\iK_{i}$ and the mesh point of $K_{i}^{E}$ (and $K_{i}^{F}$) with the midpoint of $E_{i}$ (and $F_{i}$ respectively), note that in this case $M_{i}$ has $i$ as its vertex as well. 
\end{enumerate}
The collection of all $M_{i}$ defines a dual mesh $ \iT_{h}^{*} $ such that
$$ \overline{\Omega} = \bigcup_{i \in \iN} M_{i}.
$$
As for $ \iT_{h} $, we denote by $ C_{M} $ the mesh point of $ M \in \iT_{h}^{*}$. Note that if $ M_{i} $ has edges lying on $ \partial \Omega $, then $ C_{M_{i}} $ is chosen to be the corresponding vertex $i$ of the primal mesh (see Figure). Finally, we construct a third grid $ \iT_{h}^{**} $ as a triangular subgrid of the dual grid as follows: for an element $ M \in \iT_{h}^{*} $, we construct elements of $ \iT_{h}^{**} $ by connecting $ C_{M} $ to all vertices of $ \iT_{h}^{*} $ (see Figure):
$$ \overline{\Omega} = \bigcup_{T \in \iT_{h}^{**}} \overline{T}.
$$

\begin{remark} \label{rmk:macro}
By construction, each dual control volume $M \in \iT_{h}^{*}$ is indeed a macroelement - the union of a fixed number of adjacent elements of the third mesh. If one chooses $c_{M_{i}}$ to be the associated vertex $i$ (of the primal mesh), then $M_{i}$ consists of triangles $T \in \iT_{h}^{**}$ that have $i$ as their vertex. 
\end{remark}

By applying the FECC scheme, the finite element space for the displacement $ \bu $ is the standard finite elements of order 1 defined on the third mesh $ \iT_{h}^{**} $. For the pressure $p$, using macroelement techniques \cite{RS84}, \cite[pp.235-238] {BF91} and due to Remark~\ref{rmk:macro}, the finite element space for $p$ is chosen to be P0 functions on the dual mesh. Thus, $p_{h}$ is piecewise constant on each macroelement $M \in \iT_{h}^{*}$. The effect of this choice will be show in the next section when we prove the stability of the resulting scheme. In particular:
\begin{enumerate}
\item[(a)] For the pressure:
$$ p_{h} \in Q_{h}=\left \{ q_{h} \in L_{0}^{2}(\Omega): q_{h}\vert_{M} \in \mathbb{P}^{0}(M), \; \forall M \in \iT_{h}^{*} \right \} \subset L_{0}^{2}(\Omega).
$$
Thus
\begin{equation} \label{FEp}
 p_{h}(\bx) = \sum_{M \in \iT_{h}^{*}} p_{M} \chi_{M}(\bx),
\end{equation}
where $ \chi_{M} $ is the characteristic function of $ M \in \iT_{h}^{*} $. 
\item[(b)] For the displacement: 
$$ \bu_{h} \in \bV_{h}=\left \{ \bv_{h} \in \bV_{0}: \bv_{h} \vert_{T} \in \left (\mathbb{P}^{1}(T)\right )^{2}, \; \forall T \in \iT_{h}^{**}  \right \} \subset \bV_{0}.
$$
The basis functions of $ \bV_{h} $ are defined at the nodes of elements of the third mesh. Since we impose homogeneous Dirichlet boundary conditions, we only need to deal with interior nodes. Denote by $\iN^{**}$ the set of interior nodes of elements of $\iT_{h}^{**}$.
\begin{remark} \label{rmk:nodes}
 By construction, $\iN^{**}$ consists of mesh points of the primal mesh and mesh points of interior dual control volumes:
 $$\iN^{**} = \bigcup_{K \in \iT_{h}} C_{K} \cup \bigcup_{\substack{M \in \iT_{h}^{*} \\ \partial M \cap \partial \Omega = \emptyset}} C_{M}.$$
\end{remark}
Denote by $N_{P}$ the basis function of $\bV_{h}$ at node $P \in \iN^{**}$, we seek for $ u_{h} \in V_{h} $ of the form:
\begin{align} 
\bu_{h}(\bx)&= \sum_{P \in \iN^{**}} \left (u_{P}^{(1)}N_{P}(\bx), u_{P}^{(2)}N_{P}(\bx)\right ), \label{FEu} 
\end{align}
where $ \bu_{P} = \left (u_{P}^{(1)}, u_{P}^{(2)}\right ) $ is the nodal values of $ \bu_{h} $ at the vertex $ P \in \iN^{**}$. 
\end{enumerate}
The discrete variational formulation of problem \eqref{variational-mixed} is then 
\begin{eqnarray} \label{discreteVM}
\mbox{Find $ \bu_{h} \in \bV_{h} $ and $ p_{h} \in Q_{h} $ such that}
   \hspace{4cm}\nonumber\\    
\begin{array}{rll} a(\bu_{h}, \bv_{h}) + b(\bv_{h}, p_{h})   & = L_{f}(\bv_{h} ), & \forall \bv_{h} \in \bV_{h}, \vspace{3pt}\\
 b(\bu_{h}, q_{h}) - \frac{1}{\lambda } c(p_{h},q_{h}) & = 0, & \forall q_{h} \in Q_{h}.
\end{array}
\end{eqnarray}
\subsection{Well-posedness of the discrete problem} \label{Subsec:proof}
According to the theory of mixed finite elements \cite[Chapter II, \S 2]{BF91}, the well-posedness of \eqref{discreteVM} is given by the following three conditions 
\begin{enumerate}
\item The bilinear form $a(\cdot, \cdot)$ is continuous, symmetric on $\bV_{h} \times \bV_{h}$ and is uniformly coercive on $\bV_{h}^{0}:=\left \{ \bv_{h} \in \bV_{h}: b(\bv_{h}, q_{h}) = 0, \forall q_{h} \in Q_{h} \right \}$, i.e. there exists $\alpha >0$ independent of the mesh size $h$ such that 
$$ a(\bv_{h}, \bv_{h}) \geq \alpha_{0} \| \bv_{h} \|_{1}, \; \forall \bv_{h} \in \bV_{h}^{0}. 
$$
\item The bilinear form $b(\cdot, \cdot)$ is continuous on $\bV_{h} \times Q_{h}$ and satisfies the uniform inf-sup condition (or Babu\v{s}ka-Brezzi stability condition), i.e. there exists $\beta >0$ independent of $h$ such that 
$$ \inf_{q_{h} \in Q_{h}} \, \sup_{\bv_{h} \in \bV_{h}} \frac{b(\bv_{h}, q_{h})}{\|\bv_{h}\|_{1} \| q_{h}\|_{0}} \geq \beta, \; \| \bv_{h}\|_{1} \neq 0, \, \| q_{h} \|_{0} \neq 0. 
$$
\item The bilinear form $c(\cdot, \cdot)$ is continuous, symmetric on $Q_{h} \times Q_{h}$ and is positive semi-definite:
$$  c(q_{h}, q_{h}) \geq 0, \; \forall q_{h} \in Q_{h}. 
$$ 
\end{enumerate}
We now check these three conditions: the continuity of the bilinear forms $a(\cdot, \cdot)$, $b(\cdot, \cdot)$ and $c(\cdot, \cdot)$ on their associated spaces is straightforward and so the symmetry of $a(\cdot, \cdot)$ and $c(\cdot, \cdot)$. It is simple to verify that $c(\cdot, \cdot)$ is positive semi-definite. The uniform coercivity of $a(\cdot, \cdot)$ is given by Korn's first inequality inequality \cite[Chapter III.3]{DL76}. There only remains to show that the $b(\cdot, \cdot)$ satisfies the uniform inf-sup condition. This is usually the main problem one has to deal with in order to prove the stability of the numerical scheme. However, in our case, this can be obtained directly by using the macroelement condition \cite{RS84}. In particular, for a macroelement $M \in \iT_{h}^{*}$, define the spaces
$$ \bV_{0,M}:= \left \{ \bv_{h} \in \bV_{h}: \bv_{h} = 0 \; \text{in} \; \Omega \setminus M \right \},
$$
and 
$$ N_{M}:=\left \{ q_{M}: q_{M} = q_{h} \vert_{M}, q_{h} \in Q_{h}, \int_{M} q_{h} \Div \bv_{h} d\bx = 0, \; \forall \bv_{h} \in \bV_{0,M}\right  \}. 
$$
By the definitions of $\bV_{h}$ and $Q_{h}$, we have that $N_{M}$ is one-dimensional (indeed, $q_{M}$ is piecewise constant on $M$ and $q_{M} \vert_{T_{1}} = q_{M} \vert_{T_{2}}, \; \forall T_{1}, T_{2} \subset M$).  
  
Using classical results of the approximation of the saddle point problem, we have the following theorems:
\begin{theorem} \label{thrm:well-posed}
There exists a unique solution $(\bu_{h}, p_{h}) \in \bV_{h} \times Q_{h}$ of the discrete problem \eqref{discreteVM}. 
\end{theorem}
\begin{theorem}\label{thrm:estimates}
Let $(\bu, p) \in \bV_{0} \times L_{0}^{2}(\Omega)$ and $(\bu_{h}, p_{h}) \in \bV_{h} \times Q_{h}$ be the solutions to problems \eqref{variational-mixed} and \eqref{discreteVM} respectively. Then the following estimate, which is uniform with respect to $\lambda$, hold:
$$ \| \bu - \bu_{h} \|_{1} + \| p - p_{h} \|_{0} \leq C \left ( \inf_{\bv_{h} \in \bV_{h}} \| \bu - \bv_{h} \|_{1} + \inf_{q_{h} \in Q_{h}} \| p- q_{h} \|_{0} \right ),
$$
where $C$ is a constant independent of the mesh size. 
\end{theorem}  
In the rest of this section, we will present in detail the calculation of the scheme \eqref{discreteVM} for implementation purpose. It will be shown that the scheme is cell-centered in the sense that the unknowns for $\bu$ represent the average values of $\bu$ over primal elements $K \in \iT_{h}$ (after performing some linear transformation) while the unknowns for $p$ represent the average values of $p$ over dual control volumes $M \in \iT_{h}^{*}$. Thus, the number of unknowns for $\bu$ is the number of elements of the primal mesh and that for $p$ is the number of nodes of elements of the primal mesh. However, the scheme is more effective than other cell-centered scheme (using the same order of the approximation space) because we have approximated $\bu$ on a finer mesh (the third mesh) without requiring additional computational cost.
\subsection{Accuracy and computational cost of FECC scheme}\label{Subsec:efficiency}
In order to simplify the calculation process, we firstly rewrite \eqref{discreteVM} equivalently as \cite[Chapter VI, p.201]{BF91}:
\begin{eqnarray} \label{discreteVM2}
\mbox{Find $ \bu_{h} \in \bV_{h} $ and $ p_{h} \in Q_{h} $ such that}
   \hspace{4cm}\nonumber\\    
\begin{array}{rll} \mu \int_{\Omega} \nabla \bu_{h} : \nabla \bv_{h} \, d \bx + (\mu +1) \int_{\Omega} p_{h} \Div \bv_{h}  \, d\bx  & = \int_{\Omega} \bff \cdot \bv_{h} \, d \bx, & \forall \bv_{h} \in \bV_{h}, \vspace{7pt}\\ 
 \int_{\Omega} q_{h} \Div \bu_{h} \,  d \bx - \frac{1}{\lambda } \int_{\Omega} p_{h} \, q_{h} \, d \bx & = 0, & \forall q_{h} \in Q_{h}.
\end{array}
\end{eqnarray}
In the following, we make use of the notation
$$ \partial_{x_{1}} g:= \frac{\partial g}{\partial x_{1}}, \; \partial_{x_{2}} g:= \frac{\partial g}{\partial x_{2}}, 
$$
for some function $g \in H^{1}(\Omega)$. 

Now to obtain the linear system from \eqref{discreteVM2}, we proceed as usual by choosing test functions of the first and the second equations of \eqref{discreteVM} as basis functions of $\bV_{h}$ and $Q_{h}$ respectively. 

\begin{itemize}
\item[(a)] From the representations of $\bu_{h}$ and $p_{h}$ in \eqref{FEu} and \eqref{FEp} respectively, letting $\bv_{h} = \left (N_{Q}(\bx), 0\right ), \; Q \in \iN^{**}$ in the first equation of \eqref{discreteVM2} to obtain
\begin{multline} \label{eq2:u1}
\mu \int_{\Omega} \sum_{P \in \iN^{**}}  \left ( u_{P}^{(1)} \partial_{x_{1}} N_{P}(\bx) \; \partial_{x_{1}} N_{Q}(\bx) + u_{P}^{(1)} \partial_{x_{2}} N_{P}(\bx) \; \partial_{x_{2}} N_{Q}(\bx) \right ) d\bx + \\
 (\mu +1) \int_{\Omega} \partial_{x_{1}} N_{Q}(\bx) \; \sum_{M \in \iT_{h}^{*}} p_{M} \chi_{M}(\bx)d\bx = \int_{\Omega} f^{(1)} \, N_{Q}(\bx) d\bx, \; \;  \forall Q \in \iN^{**}.
\end{multline}
Denote by $\iN_{M}$ the set of nodes of elements $T \in \iT_{h}^{**}$ that $T \subset M$. 
In the following, we will show how to condense out the unknowns $u_{C_{M}}^{(1)}, \; \forall M~\in~\iT_{h}^{*}, \; \partial M \cap \partial \Omega = \emptyset$ (i.e. interior dual control volumes) from the formulation. In fact, $u_{C_{M}}^{(1)}$ can be computed by a linear combination of unknowns at primal mesh points $u_{C_{K}}^{(1)}, \, K \in \iT_{h} $. \vspace{4pt}\\ To this purpose, let $\bv_{h} = \left (N_{C_{M}}, 0 \right ), \; M \in \iT_{h}^{*}, \; \partial M \cap \partial \Omega = \emptyset$ then \eqref{eq2:u1} becomes
\begin{multline} \label{eq2:u1bis}
\mu \int_{M} \sum_{P \in \iN_{M}}  \left ( u_{P}^{(1)} \partial_{x_{1}} N_{P}(\bx) \; \partial_{x_{1}} N_{C_{M}}(\bx) +  u_{P}^{(1)} \partial_{x_{2}} N_{P}(\bx) \; \partial_{x_{2}} N_{C_{M}}(\bx) \right ) + \\
(\mu +1) \int_{M} \partial_{x_{1}} N_{C_{M}}(\bx) \;  p_{M} = \int_{M} f^{(1)} \, N_{C_{M}}(\bx) d\bx.
\end{multline}
According to Remark~\ref{rmk:nodes}, $\iN_{M}$ consists of the mesh point $C_{M}$ of $M$ and the set of mesh points $C_{K}$ of elements $K \in \iT_{h}$ such that $K~\cap~M~\neq~\emptyset$.  
Thus, from \eqref{eq2:u1bis} we can compute $u_{C_{M}}^{(1)}$ using unknowns at primal mesh points of elements that intersect $M$:
\begin{multline} \label{eq2:uM1}
u_{C_{M}}^{(1)} = \Pi_{M}^{(1)} \left ( \left \{ u_{C_{K}}\right \}_{K \in \iT_{h}, \, K \cap M \neq \emptyset}, \bff \right )  =\\
 -\frac{1}{\Theta_{M}} \int_{M} \sum_{\substack{K \in \iK_{h} \\ K \cap M \neq \emptyset}}  \left ( u_{C_{K}}^{(1)} \partial_{x_{1}} N_{C_{K}}(\bx) \; \partial_{x_{1}} N_{C_{M}}(\bx) +  u_{C_{K}}^{(1)} \partial_{x_{2}} N_{C_{K}}(\bx) \; \partial_{x_{2}} N_{C_{M}}(\bx) \right ) + \\
 \frac{1}{\mu \, \Theta_{M}} \left ((\mu + 1) \int_{M} \partial_{x_{1}} N_{C_{M}}(\bx) \;  p_{M} - \int_{M} f^{(1)} \, N_{C_{M}}(\bx) d\bx \right ),
\end{multline}
where 
$$ \Theta_{M}:= \int_{M} \left (\partial_{x_{1}} N_{C_{M}}\right )^{2} + \left (\partial_{x_{2}} N_{C_{M}}\right )^{2}  d \bx. 
$$
For convenience, we write  $\Pi_{M}^{(1)}$ for 
 $\Pi_{M}^{(1)} \left ( \left \{ u_{C_{K}}\right \}_{K \in \iT_{h}, \, K \cap M \neq \emptyset}, \bff \right )$.
%
Now substituting \eqref{eq2:uM1} into \eqref{eq2:u1} we obtain
\begin{multline} \label{eq2:u1Final}
\mu \int_{\Omega} \sum_{K \in \iT_{h}}  \left ( u_{C_{K}}^{(1)} \partial_{x_{1}} N_{C_{K}}(\bx) \; \partial_{x_{1}} N_{C_{H}}(\bx) + u_{C_{K}}^{(1)} \partial_{x_{2}} N_{C_{K}}(\bx) \; \partial_{x_{2}} N_{C_{H}}(\bx) \right ) d\bx + \\
\mu \int_{\Omega} \sum_{\substack{M \in \iT_{h}^{*} \\ \partial M \cap \partial \Omega = \emptyset}}  \left ( \Pi_{M}^{(1)} \partial_{x_{1}} N_{C_{M}}(\bx) \; \partial_{x_{1}} N_{C_{H}}(\bx) + \Pi_{M}^{(1)} \partial_{x_{2}} N_{C_{M}}(\bx) \; \partial_{x_{2}} N_{C_{H}}(\bx) \right ) d\bx + \\
(\mu + 1)  \int_{\Omega} \partial_{x_{1}} N_{C_{H}}(\bx) \; \sum_{M \in \iT_{h}^{*}} p_{M} \chi_{M}(\bx)d\bx = \int_{\Omega} f^{(1)} \, N_{C_{H}}(\bx) d\bx, \; \;  \forall H \in \iT_{h}.
\end{multline}
This linear equation only involves unknowns of the primal elements. 
\item[(b)] Similar to (a), we now choose the test function $\bv_{h} = \left (0, N_{Q}(\bx) \right ), \; Q \in \iN^{**}$ in the first equation of \eqref{discreteVM} and obtain
\begin{multline} \label{eq2:u2Final}
\mu \int_{\Omega} \sum_{K \in \iT_{h}}  \left ( u_{C_{K}}^{(2)} \partial_{x_{2}} N_{C_{K}}(\bx) \; \partial_{x_{2}} N_{C_{H}}(\bx) + u_{C_{K}}^{(2)} \partial_{x_{1}} N_{C_{K}}(\bx) \; \partial_{x_{1}} N_{C_{H}}(\bx) \right ) d\bx + \\
\mu \int_{\Omega} \sum_{\substack{M \in \iT_{h}^{*} \\ \partial M \cap \partial \Omega = \emptyset}}  \left ( \Pi_{M}^{(2)} \partial_{x_{2}} N_{C_{M}}(\bx) \; \partial_{x_{2}} N_{C_{H}}(\bx) + \Pi_{M}^{(2)} \partial_{x_{1}} N_{C_{M}}(\bx) \; \partial_{x_{1}} N_{C_{H}}(\bx) \right ) d\bx + \\
 \int_{\Omega} \partial_{x_{2}} N_{C_{H}}(\bx) \; \sum_{M \in \iT_{h}^{*}} p_{M} \chi_{M}(\bx)d\bx = \int_{\Omega} f^{(2)} \, N_{C_{H}}(\bx) d\bx, \; \;  \forall H \in \iT_{h},
\end{multline}
where 
\begin{multline*}
 \Pi_{M}^{(2)} = \Pi_{M}^{(2)} \left ( \left \{ u_{C_{K}}\right \}_{K \in \iT_{h}, \, K \cap M \neq \emptyset}, \bff \right ) =\\
 -\frac{1}{\Theta_{M}} \int_{M} \sum_{\substack{K \in \iK_{h} \\ K \cap M \neq \emptyset}}  \left ( u_{C_{K}}^{(2)} \partial_{x_{1}} N_{C_{K}}(\bx) \; \partial_{x_{1}} N_{C_{M}}(\bx) +  u_{C_{K}}^{(2)} \partial_{x_{2}} N_{C_{K}}(\bx) \; \partial_{x_{2}} N_{C_{M}}(\bx) \right ) + \\
 \frac{1}{\mu \, \Theta_{M}} \left ((\mu + 1) \int_{M} \partial_{x_{2}} N_{C_{M}}(\bx) \;  p_{M} - \int_{M} f^{(2)} \, N_{C_{M}}(\bx) d\bx \right ).
\end{multline*}
\item[(c)] Letting $q_{h} = q_{M}, \, M \in \iT_{h}^{*}$, we rewrite the second equation of \eqref{discreteVM} as: 
\begin{equation*}
\sum_{P \in \iN_{M}} \int_{M} u_{P}^{(1)} \partial_{x_{1}} N_{P}(\bx) + u_{P}^{(2)} \partial_{x_{2}} N_{P}(\bx) d \bx -\frac{1}{\lambda} \int_{M} p_{M} d\bx = 0, \; \forall M \in \iT_{h}^{*}.
\end{equation*}
Applying the results in $(a)$ and $(b)$ we have
\begin{multline} \label{eq1:pFinal}
\sum_{\substack{K \in \iT_{h} \\ K \cap M \neq \emptyset}} \int_{M} u_{C_{K}}^{(1)} \partial_{x_{1}} N_{C_{K}}(\bx) + u_{C_{K}}^{(2)} \partial_{x_{2}} N_{C_{K}}(\bx) d \bx  + \\
\int_{M} \Pi_{M}^{(1)} \partial_{x_{1}} N_{C_{M}}(\bx) + \Pi_{M}^{(2)} \partial_{x_{2}} N_{C_{M}}(\bx) d \bx -\frac{1}{\lambda} \int_{M} p_{M} d\bx = 0, \; \forall M \in \iT_{h}^{*}.
\end{multline}
\end{itemize}
Finally, we end up with a linear system of problem \eqref{discreteVM} as follows:
\begin{equation}
\left [ \begin{array}{ll} \bA & \bB^{t} \\
\bB & -\frac{1}{\lambda} \bC
\end{array} \right ] \, \left [ \begin{array}{c} \bU \\
\bP
\end{array} \right ] =  \left [ \begin{array}{c} \bF \\
0
\end{array} \right ],
\end{equation}
where $\bU$ is a vector $\left ( u_{C_{T}}^{(1)}, u_{C_{T}}^{(2)} \right )_{T \in \iT_{h}}$ and $\bP$ is $\left (p_{M}\right )_{M \in \iT_{h}^{*}}$, which implies that the scheme is cell-centered.

%
%
\section{Numerical experiments}
\label{Sec:NumRe}
We carry out several test cases in 2D to verify the performance of the constructed scheme and to compare its performance with other schemes for nearly incompressible elasticity problems. 

\section*{References}



\end{document}